\documentclass[11pt]{article}

\usepackage[margin=1.1in]{geometry}
\usepackage{graphicx}
\usepackage{tikz}
\usetikzlibrary{arrows.meta}

\usepackage{amssymb}

\newcommand{\bfa}{{\bf a}}
\newcommand{\bfb}{{\bf b}}
\newcommand{\bfc}{{\bf c}}

\newtheorem{theorem}{Theorem}[section]

\newtheorem{lemma}[theorem]{Lemma}

\begin{document}

\title{Extrapolating on Taylor Series Solutions \\
       of Homotopies with Nearby Poles\thanks{Supported by 
the National Science Foundation under grant DMS 1854513.}}

\author{Jan Verschelde\thanks{University of Illinois at Chicago,
Department of Mathematics, Statistics, and Computer Science,
851 S. Morgan St. (m/c 249), Chicago, IL 60607-7045.
Email: {\tt janv@uic.edu}, URL: {\tt http://www.math.uic.edu/$\sim$jan}.}
\and Kylash Viswanathan\thanks{University of Illinois at Chicago,
Department of Mathematics, Statistics, and Computer Science,
851 S. Morgan St. (m/c 249), Chicago, IL 60607-7045.
Email: {\tt kviswa5@uic.edu}.}}

\date{26 April 2024}

\maketitle

\begin{abstract}
A polynomial homotopy is a family of polynomial systems in one parameter,
which defines solution paths starting from known solutions and ending
at solutions of a system that has to be solved.
We consider paths leading to isolated singular solutions,
to which the Taylor series converges logarithmically.
Whether or not extrapolation algorithms manage to accelerate
the slowly converging series depends on the proximity of poles
close to the disk of convergence of the Taylor series.

{\bf Keywords.} continuation, extrapolation, homotopy, pole,
polynomial, Taylor, series, singularity.
\end{abstract}
\section{Introduction}

Polynomial homotopies define solution paths leading to solutions
of polynomial systems.  Many interesting systems have singular solutions
and then we often observe the convergence of Taylor series expansions 
to be logarithmic.
Logarithmic convergence is extremely slow, as obtaining one extra bit
of accuracy costs a multiple of the cost of the current accurate bits.

Unless certain conditions hold, sequences that converge logarithmically
cannot be accelerated~\cite{DGB80,DGB82}, \cite{Kow81}.
We illustrate the construction of Taylor series for which extrapolation
algorithms fail to accelerate the sequence.  In particular, we show
that the proximity of another nearby pole, near the disk of convergence
of the Taylor series, determines whether Richardson extrapolation~\cite{Bre80},
Aitken extrapolation~\cite{BZ91}, \cite{Sid03}, \cite{Wen89},
the rho algorithm~\cite{Wyn56} and theta algorithm~\cite{Bre71}
are effective.

This paper follows the prior work in~\cite{BV18}, \cite{TVV20},
and~\cite{VV22}.  The results in~\cite{VV22} stated the error
expansions in~$1/n$ as necessary for the effectiveness of 
Richardson extrapolation.  We consider special situations
where error expansions in~$1/n$ cannot be constructed.
Computational experiments for this paper were done mainly with phcpy,
the Python interface to PHCpack~\cite{Ver99},~\cite{VV24}.

\section{Problem Statement}

Our problem originates from a calculation on an example 
for which the rho algorithm performed spectacularly well, 
followed then by the complete failure on a slightly more
complicated example.

\subsection{A Motivating Example}

Consider the solutions $x(t)$ of
\begin{equation}
   x^q = (1 - t)^p, \quad t \in [0, 1], \quad p \geq 2, \quad q \geq 1.
\end{equation}
At $t=0$, one solution path starts at $x = 1$ and ends 
at $t=1$, at a singularity.

The Taylor series of $x(t) = a_0 + a_1 t + a_2 t^2 + \cdots$ 
has nonzero coefficients $a_n$ for all $n$.  So, the ratio
$a_n/a_{n+1}$ is well defined and converges very slowly to~1.
To accelerate the convergence of $a_n/a_{n+1}$,
we apply the rho algorithm, as defined in the Python function below.

\begin{verbatim}
def rhoAlgorithm(nbr):
    """
    Runs the rho algorithm in rational arithmetic,
    on the numbers given in the list nbr, using x(n) = n+1.
    Returns the last element of the table of extrapolated numbers.
    """
    rho1 = [1/(nbr[n] - nbr[n-1]) for n in range(1, len(nbr))]
    rho = [nbr, rho1]
    for k in range(2, len(nbr)):
        nextrho = []
        for n in range(k, len(nbr)):
            invrho1 = k/(rho[k-1][n-k+1] - rho[k-1][n-k])
            nextrho.append(rho[k-2][n-k+1] + invrho1)
        rho.append(nextrho)
    return nextrho[-1]
\end{verbatim}
For example, for $p=7$ and $q = 18$, the first three numbers in 
the sequence $a_n/a_{n+1}$ are $36/11$, $54/29$, and $72/47$.
It turns out that this input suffices for the rho algorithm to
return~1.  A systematic run for all integers $q$ ranging between 2 and 20,
and for all integers $p$ from~1 to $q-1$ returns~1 as well.

One possible justification of the good performance of the rho algorithm
on this example is the connection with Thiele interpolation~\cite{CW87},
as the rho algorithm is designed to compute even order convergents
of Thiele's interpolating continued fraction~\cite{Wen89}.
Observe the inverse divided differences in the Python function above,
which could cause division by zero and other numerical errors.
The rules given in~\cite{Red92} are shown to improve the numerical
stability of extrapolation algorithms.

In this motivating example, only four terms in a Taylor series suffice 
for the rho algorithm to accurately compute the pole, 
which happens to be the only pole.
The next section provides a slightly more complicated situation.

\subsection{Homotopies with Prescribed Poles}

\begin{figure}[hbt]
\begin{center}
\begin{tikzpicture}
\draw (0cm, 0cm) circle(1cm);
\draw (0mm, -2mm) -- (0mm, 2mm);
\draw (0mm, 0mm) -- (-5mm, 10mm);
\draw[fill=red] (1cm, 0cm) circle(0.7mm);
\draw[-{Latex[length=2mm, width=1mm]}] (-15mm, 0mm) -- (17mm, 0mm);
\draw[fill=orange] (-5mm, 10mm) circle(0.7mm);
\node[text width=5mm] at (-6mm, 11.34mm) {$P$};
\node[text width=5mm] at (20mm, -1mm) {$t$};
\node[text width=5mm] at (13mm, 2mm) {$1$};
\node[text width=5mm] at (1.5mm, -4mm) {$0$};
\draw[fill=green] (0mm, 0mm) circle(0.7mm);
\end{tikzpicture}
\end{center}
\caption{One nearby pole $P$ near the disk with radius~$1$.}
\label{fignearbypoles}
\end{figure}
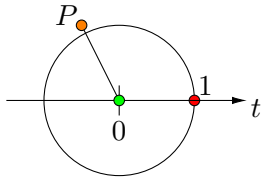

Figure~\ref{fignearbypoles} corresponds to the monomial homotopy
for $P = -1/2 + I$, $I = \sqrt{-1}$,
\begin{equation} \label{eqnearbypoles}
  x^2 - C (1 - t) (P - t) = 0, \quad t \in [0, 1],
\end{equation}
where the constant $C$ is such that $x(0) = \pm 1$.
There is a double root at $t = 1$ and at $P$.
$P$ is near the disk of convergence of the Taylor series of $x(t)$.

Although the homotopy in~(\ref{eqnearbypoles}) is only slightly
more complicated than the motivating example,
extrapolation methods are no longer effective.
It is the smallest illustration for our problem statement: 
how close should $P$ be to the disk of convergence
for extrapolation to be effective?

\subsection{Fabry's Ratio Theorem and Extrapolation}

All Taylor series considered in this paper have nonzero coefficients
and the ratio theorem of Fabry can be applied.

\begin{theorem} 
[\mbox{\rm the ratio theorem of Fabry~\cite{Fab1896}}] \label{thmFabry}

\noindent {\em If for the series
\begin{equation}
   x(t) = c_0 + c_1 t + c_2 t^2 + \cdots + c_n t^n + c_{n+1} t^{n+1} + \cdots,
\end{equation}
we have
\begin{equation}
   \lim_{n \rightarrow \infty} \frac{c_n}{c_{n+1}} = P,
\end{equation}
then
\begin{itemize}
\item $P$ is a singular point of the series, 
\item it lies on the boundary of the circle of convergence of the series, and
\end{itemize}
then the radius of the disk of convergence is~$|P|$.}
\end{theorem}
In case the assumptions of Theorem~\ref{thmFabry} hold,
then with a coordinate transformation the nearest singularity can 
be placed at $t=1$.  If the radius of convergence equals one,
then $|c_{n+1}| \approx |c_n|$, for sufficiently large values of~$n$.
For radius $R > 1$, the coefficients of the series are decaying
and we have $|c_{n+1}| \approx |c_n|/R$.

For Richardson extrapolation~\cite{BZ91},~\cite{Sid03} to be effective,
in~\cite{VV22}, the following was derived:
\begin{equation} \label{eqderived}
  \left| 1 - \frac{c_n}{c_{n+1}} \right|
  = \sum_{k=1}^\infty \gamma_k \left( \frac{1}{n} \right)^k,
\end{equation}
for complex constants $\gamma_k$.
This derivation was based on experiments on the Taylor series
of $\sqrt{1-t}$.  Extrapolation methods were observed not to work
when there was nearby pole, around the same radius one.
But when the other pole was farther away, with a radius~$R$
of a larger magnitude than one, e.g.~$R > 10$, then extrapolation
methods were again effective.

Experiments have shown that, if Richardson extrapolation works,
then Aitken extrapolation works too, 
and the rho and theta algorithms work as well, often much better.
Therefore, we focus on series expansions in $1/n$,
which provides the justification for the effectiveness
of Richardson extrapolation.

\section{Error Expansions for Nearby Poles}

We start by formalizing the motivating example in the first lemma.
\begin{lemma} \label{lemma0}
Let $a_n$ by the $n$-th coefficient of the Taylor series $x(t)$ defined by
\begin{equation}
   x^q = (1 - t)^p, \quad t \in [0, 1], \quad p \geq 2, \quad q \geq 1.
\end{equation}
Then:
\begin{equation} \label{eqlemma0result}
   \frac{a_n}{a_{n+1}}
   = 1 + \sum_{k=1}^\infty \alpha_k \left( \frac{1}{n} \right)^k,
\end{equation}
for coefficients $\alpha_k$ independent of~$n$.
\end{lemma}
{\em Proof.}
The Taylor series of $x(t)$ is
\begin{equation}
   x(t) = x(0) + x'(0) + \frac{x''(0)}{2} + \cdots
        + \frac{x^{(n)}(0)}{n!} + \cdots.
\end{equation}
The $n$-th derivative of $x(t) = (1 - t)^{p/q}$ is
\begin{equation}
   x^{(n)}(t) = \frac{p}{q} 
   \left( \frac{p}{q} - 1 \right) \cdots
   \left( \frac{p}{q} - n + 1 \right) (1 - t)^{p/q - n} (-1)^n.
\end{equation}
The ratio of two consecutive Taylor series coefficients is then
\begin{eqnarray}
   \frac{a_n}{a_{n+1}} 
   & = & {\displaystyle
          \frac{\frac{x^{(n)}(0)}{n!}}
               {\frac{x^{(n+1)}(0)}{(n+1)!}}
   ~ = ~  \frac{(n+1)!}{n!}
          \frac{x^{(n)}(0)}{x^{(n+1)}(0)}} \\
   & = & {\displaystyle 
          \left( n+1 \vphantom{\frac{1}{2}} \right)
          \left( \frac{-1}{\frac{p}{q} - n} \right)
   ~ = ~ \left( n+1 \vphantom{\frac{1}{2}} \right)
       \left( \frac{1}{ n 
       \left( 1 - \frac{p}{q} \left( \frac{1}{n} \right) \right)} \right)} \\
   & = & {\displaystyle \frac{n+1}{n} 
       \left( \frac{1}{ 1 - \frac{p}{q} \left( \frac{1}{n} \right)} \right)
   ~ = ~ \left( 1 + \frac{1}{n} \right)
       \left( \frac{1}{ 1 - \frac{p}{q} \left( \frac{1}{n} \right)} \right).}
\end{eqnarray}
Applying the geometric series expansion to the last factor and 
collecting the terms in~$1/n$ gives~(\ref{eqlemma0result}). \hfill Q.E.D.

A Taylor series development is unique and Taylor series are invertible.
If we have a series with two nearby poles, then we can view the series as
a product of two Taylor series, each series having only one nearby pole.
The second lemma below rewrites the ratio of the coefficients of the product
in terms of ratios of the coefficients of the factors in the product.

\begin{lemma} \label{lemma1} Let $\bfa = (a_k)_{k=0}^\infty$ 
and $\bfb = (b_k)_{k=0}^\infty$ be sequences of nonzero numbers.
Assume the convolution $\bfc = \bfa \star \bfb$ yields another sequence
of nonzero numbers $\bfc = (c_k)_{k=0}^\infty$.
Then we have
\begin{equation} \label{eqlemma1result}
   \frac{c_n}{c_{n+1}} =
   \left( \frac{a_n}{a_{n+1}} \right)
   \sum_{k=0}^n \frac{1}{\displaystyle \sum_{\ell=0}^{n+1} 
      \left( \frac{a_n}{a_{n+1}} \right)
      \left( \frac{a_\ell}{a_k} \right)
      \left( \frac{b_{n+1-\ell}}{b_{n-k}} \right)}.
\end{equation}
\end{lemma}
{\em Proof.}
Following the definition of the convolution
\begin{eqnarray}
   \frac{c_n}{c_{n+1}}
   & = & \frac{\displaystyle \sum_{k=0}^n a_k b_{n-k}}
          {\displaystyle \sum_{\ell=0}^{n+1} a_\ell b_{n+1-\ell}}
   ~ = ~ \sum_{k=0}^n \frac{a_k b_{n-k}}
          {\displaystyle \sum_{\ell=0}^{n+1} a_\ell b_{n+1-\ell}} \\
   & = & \sum_{k=0}^n \frac{1}
          {\frac{\displaystyle
             \sum_{\ell=0}^{n+1} a_\ell b_{n+1-\ell}}{a_k b_{n-k}}}
   ~ = ~ \sum_{k=0}^n \frac{1}
          {\displaystyle
             \sum_{\ell=0}^{n+1} \frac{a_\ell b_{n+1-\ell}}{a_k b_{n-k}}}.
\end{eqnarray}
In the denominator we then extract the factor $a_{n+1}/a_n$:
\begin{equation}
   \frac{c_n}{c_{n+1}}
   = \sum_{k=0}^n \frac{1}
          {\displaystyle
             \left( \frac{a_{n+1}}{a_n} \right)
             \sum_{\ell=0}^{n+1} 
             \left( \frac{a_n}{a_{n+1}} \right)
             \left( \frac{a_\ell}{a_k} \right)
             \left( \frac{b_{n+1-\ell}}{b_{n-k}} \right)}
\end{equation}
which then yields~(\ref{eqlemma1result}). \hfill Q.E.D.

The result of Lemma~\ref{lemma1} hints at the relative importance
of the coefficients in the developments of $c_n/c_{n+1}$.
If on the one hand, the ratios of the first factor~$\bfa$ dominate,
then we may expect that, if extrapolation works on those first ratios,
then it will also work on the product.
If on the other hand, the ratios of both factors $\bfa$ and $\bfb$
are equally important, then the development of the product lacks the
regularity necessary for the effectiveness of any extrapolation algorithm,
because the summation over $n+1$ terms of numbers of magnitude about one 
muddle up the developments in~$1/n$.

Lemma~\ref{lemma2} introduces a large pole, sufficiently far away
in magnitude from~1, so that a development in~$1/P$ can be constructed.

\begin{lemma} \label{lemma2} For $P$ such that $|P| \gg 1$, 
for indices $k > 0$ and $\ell > 0$, 
and for nonzero coefficients $a_k, a_{k-1}, \ldots, a_{-\ell+1}, a_{-\ell}$,
consider
\begin{equation}
   Q = \frac{1}{a_k P^k + a_{k-1} P^{k-1} + \cdots
              + a_{-\ell+1} P^{-\ell + 1} + a_{-\ell} P^{-\ell}.}
\end{equation}
There exists coefficients $\beta_i$ so $Q$ can be written as
\begin{equation} \label{eqlemma2result}
   Q = \sum_{i=k}^\infty \beta_i \left( \frac{1}{P} \right)^i.
\end{equation}
\end{lemma}
{\em Proof.}
We isolate the factor $a_k P^k$ in the denominator of $Q$:
\begin{equation}
   Q = \frac{1}{a_k P^k} \left(
       \frac{1}{1 + \frac{a_{k-1}}{a_k} \frac{1}{P} + \cdots
              + \frac{a_{-\ell+1}}{a_k} \frac{1}{P^{k+\ell - 1}}
              + \frac{a_{-\ell}}{a_k} \frac{1}{P^{k+\ell}}} \right).
\end{equation}
Let
\begin{equation}
   x =  \frac{a_{k-1}}{a_k} \frac{1}{P} + \cdots
      + \frac{a_{-\ell+1}}{a_k} \frac{1}{P^{k+\ell - 1}}
      + \frac{a_{-\ell}}{a_k} \frac{1}{P^{k+\ell}}
\end{equation}
and apply the geometric series
\begin{equation}
   Q = \frac{1}{a_k P^k} \left( \frac{1}{1 + x} \right) \\
     = \frac{1}{a_k P^k} \left( \sum_{i=0}^\infty (-1)^i x^i \right).
\end{equation}
Expanding $x^i$ as sums of powers of $1/P$ and collecting the terms
as the coefficients $\beta_i$ 
gives the result~(\ref{eqlemma2result}). \hfill Q.E.D.

The theorem below provides the existence of a range of~$n$
so that the development in Lemma~\ref{lemma2} can be written
as a series in~$1/n$.

\begin{theorem}  Consider
\begin{equation}
   x(t) = \sum_{k=0}^\infty a_k t^k, \quad 
   \lim_{n \rightarrow \infty} \frac{a_n}{a_{n+1}} = 1,
   \hphantom{, \quad |P| \gg 1}
\end{equation}
and
\begin{equation}
   y(t) = \sum_{k=0}^\infty b_k t^k, \quad 
   \lim_{n \rightarrow \infty} \frac{b_n}{b_{n+1}} = P, \quad 
   |P| \gg 1.
\end{equation}
Let $z(t) = x(t) \star y(t)$, with
$\displaystyle z(t) = \sum_{k=0}^\infty c_k t^k$.
Then, for sufficiently large $N$ and $M$, $M > N$, 
we have for all $n \geq N$ and $n \leq M$:
\begin{equation} \label{eqtheoremresult}
  \frac{c_n}{c_{n+1}}
   = 1 + \sum_{k=1}^\infty \gamma_k \left( \frac{1}{n} \right)^k,
\end{equation}
for constants $\gamma_k$ computed for~$n$ in $[N, M]$.
\end{theorem}
{\em Proof.} 
Applying Lemma~\ref{lemma0}, the ratio $a_n/a_{n+1}$
has an expansion in $1/n$:
\begin{equation}
   \frac{a_n}{a_{n+1}}
   = 1 + \sum_{k=1}^\infty \alpha_k \left( \frac{1}{n} \right)^k.
\end{equation}
We apply Lemma~\ref{lemma1} to $c_n/c_{n+1}$ and make the following
observations about the growth of the coefficients:
\begin{enumerate}
\item As $\displaystyle \lim_{n \rightarrow \infty} \frac{a_n}{a_{n+1}} = 1$,
      we have $a_{n+1} \approx a_n$, and for $k \geq N$ and $\ell \geq N$:
      $a_\ell \approx a_k$.
\item As $\displaystyle \lim_{n \rightarrow \infty} \frac{b_n}{b_{n+1}} = P$,
      we have $b_{n+1} \approx b_n/P \approx b_{n-1}/P^2$ and
      furthermore: $b_n \approx b_m/P^{n-m}$, for $m \geq N$.
\end{enumerate}
With these two observations, we are in position to apply Lemma~\ref{lemma2}.
Then, for $n \geq N$ and $n \leq M$, there exists some constant $C$:
\begin{equation} \label{eqmainresult}
   \frac{1}{P} = C \frac{1}{n}.
\end{equation}
The application of Lemma~\ref{lemma1} and Lemma~\ref{lemma2} gives
\begin{equation}
   \frac{c_n}{c_{n+1}} =
   \left( 1 + \sum_{k=1}^\infty \alpha_k \left( \frac{1}{n} \right)^k \right)
   \left( \sum_{i=1}^\infty \beta_i \left( \frac{1}{n} \right)^i \right)
\end{equation}
which yields~(\ref{eqtheoremresult}).  \hfill Q.E.D.

In practical examples, the sufficiently large value for~$N$
can be already as little as~8, with more typical values
about~32 or~64.  While~(\ref{eqmainresult}) suggests $|P| \approx N$,
the larger $|P|$, the better.

\section{Computational Experiments}

The rho algorithm is applied to the sequence $c_n/c_{n+1}$, 
$n=0, 1, \ldots d$,
\begin{itemize}
\item $c_n$ is the $n$-th coefficient of the Taylor series of~$x(t)$, and
\item $\displaystyle x(t) = \sqrt{\alpha~\!(1 - t) (P - t)}$, 
      where $\alpha$ is such that $x(0) = 1$.
\end{itemize}

The smallest error of the rho table for various $P$ and $d$ values
is listed in Table~\ref{taberrors}.
All calculations happened in double precision.
The coefficients $c_n$ were computed with tolerance {\tt 1.0e-12}.

\begin{table}[hbt]
\label{taberrors}
\begin{center}
\begin{tabular}{c||c|c|c}
  $P$ & $d = 8$ & $d = 16$ & $d = 32$ \\ \hline \hline
$\displaystyle -1/2 + \hphantom{16}I$  & {\tt 5.0e-01} & {\tt 3.5e-01} & {\tt 1.4e-01} \\
$\displaystyle -1/2 + \hphantom{1}2 I$ & {\tt 1.7e-01} & {\tt 9.8e-03} & {\tt 2.6e-05} \\
$\displaystyle -1\hphantom{/2} + \hphantom{1}4 I$ & {\tt 2.5e-02} & {\tt 6.9e-05} & {\tt 6.3e-09} \\
$\displaystyle -2\hphantom{/4} + \hphantom{1}8 I$ & {\tt 3.3e-03} & {\tt 5.3e-07} & {\tt 3.5e-11} \\
$\displaystyle -4\hphantom{/4} + 16 I$ & {\tt 4.1e-04} & {\tt 4.0e-09} & {\tt 2.4e-12}
\end{tabular}
\caption{Errors of the results of the rho algorithm on sequences
of size $d+1$.}
\end{center}
\end{table}

Consistent with the theoretical results, we observe that the farther
$P$ is from~1, the better the extrapolation works.

\section{Conclusions}

We provide a condition for the effectiveness of extrapolation methods 
to accelerate slowly converging Taylor series of solution paths 
of polynomial homotopies.

\noindent {\bf Acknowledgements.} 
Some of the results in this paper were presented by the first author
on 6 January 2023 in a preliminary report at the special session 
on Complexity and Topology in Computational Algebraic Geometry,
at the Joint Mathematics Meetings.
The authors thank the organizers of the session,
Ali Mohammad Nezhad and Saugata Basu, for their invitation.
The first author thanks Claude Brezinski and Michela Redivo Zaglia
for an enlighting converation at CAM23, and for pointing at~\cite{Red92}.

\bibliographystyle{plain}

\end{document}